\documentclass[12pt,a4paper]{article}

\usepackage{amsmath,amsfonts}
\usepackage{mathrsfs}

\newtheorem{theorem}{Theorem}

\def \Liminf{\mathop{\underline{\lim}}\limits}

\def\Pb{\mathbf{P}}

\def\Ex{\mathbf{E}}
\def\RR{\mathbb{R}}

\def\1{\mbox{1\hspace{-.25em}I}}
\begin{document}

\date{}

\title{On Approximation of the Backward Stochastic Differential Equation}
\author{Yury A. Kutoyants, \ \ Li Zhou\\
{\small {Laboratoire de Statistique et Processus, Universit\'e du Maine}}.\\
{\small {72085 Le Mans Cedex, France} }}

\maketitle

\begin{abstract}
We consider the problem of approximation of the solution of the backward
stochastic differential equation in the Markovian case.  We suppose that the
trend coefficient of the diffusion process depends on some unknown
 parameter and the diffusion coefficient of this equation is
small. We propose an approximation of this solution based on the
one-step MLE of the unknown parameter and we show that this
approximation is asymptotically efficient in the asymptotics of
``small noise''.
\end{abstract}


\section{Introduction}

We consider the following problem. We are given a stochastic differential
equation (called {\it forward})
\begin{equation}
\label{1}
 {\rm d} X_t=b(t,X_t)\;{\rm d} t+ a (t,X_t)\;{\rm d} W_t,\ \ X_0=x_0,\
0\leq t\leq T,
\end{equation}
and two functions $f\left(t,x,y,z\right)$  and $\Phi
\left(x\right)$. We have  to construct a couple of processes
$\left(Y_t,Z_t\right)$ such that the solution of the equation
\begin{equation}
\label{2}
   {\rm d} Y_t=-f(t,X_t,Y_t,Z_t)\;{\rm d} t+Z_t\;{\rm d} W_t,\  0\leq t\leq T,
\end{equation}
(called {\it backward}) has the final value $Y_T=\Phi \left(X_T\right)$.

The existence and uniqueness of the solution of {\it backward stochastic
  differential equation} (BSDE) in essentially more general situations was
studied by Pardoux and Peng \cite{PP1990}.  The problem \eqref{1}-\eqref{2}
considered here was introduced as forward-backward stochastic differential
equations (FBSDE) in El Karoui \& al. \cite{EPQ1997}. The solution of this
FBSDEs is presented as a triple-process $(X_t,\,Y_t,\,Z_t)_{t\geq0}$. It is
shown that the solution $(X_t,\,Y_t,\,Z_t)_{t\geq0}$ exists and is unique
under the condition that all coefficients are Lipschitzian and so on (see
\cite{EPQ1997} for details). The solution of the problem \eqref{1}-\eqref{2}
proposed in \cite{EPQ1997} is the following. Suppose that $u\left(t,x\right)$
satisfies the  equation
\begin{equation}
\label{3}
  \frac{\partial u}{\partial t}+b\left(t,x\right)\frac{\partial u}{\partial
    x}+\frac{1}{2} a \left(t,x\right)^2\frac{\partial^2 u}{\partial
    x^2}=-f\left(t,x,u, a \left(t,x\right)\frac{\partial u}{\partial x}\right)
  ,\quad u\left(T,x\right)=\Phi \left(x\right)
\end{equation}
and put  $Y_t=u\left(t,X_t\right),Z_t=a \left(t,X_t\right)
u_x'\left(t,X_t\right)$. Then by It\^o's formula the process
$Y_t$ has the stochastic differential
\begin{align*}
{{\rm d}}Y_t&=\left[\frac{\partial u}{\partial
    t}\left(t,X_t\right)+b\left(t,X_t\right)\frac{\partial u}{\partial
    x}\left(t,X_t\right)+\frac{1}{2} a
  \left(t,x\right)^2\frac{\partial^2 u}{\partial x^2}\left(t,X_t\right)
  \right]\,{{\rm d}}t\\
&\qquad \qquad + a\left(t,X_t\right)\frac{\partial u}{\partial x}\left(t,X_t\right)\,{{\rm d}}W_t\\
&=-f\left(t,X_t,Y_t,Z_t\right)\,{{\rm d}}t+Z_t\,{{\rm d}}W_t,\qquad
Y_0=u\left(0,X_0\right).
\end{align*}
The final value $Y_T=u\left(T,X_T\right)=\Phi \left(X_T\right)$. Therefore the
problem is solved and the couple $\left(Y_t,Z_t\right)$ provides the
desired solution. More details  can be found, e.g.,  in El
Karoui \& Mazliak \cite{EM} and Ma \& Yong \cite{MY}.

In the present work we consider the similar statement but in the situation
when the trend coefficient $b\left(t,x\right)$ of the diffusion process
\eqref{1} depends on the unknown parameter $\vartheta \in \Theta $, i.e.,
$b\left(t,x\right)=S\left(\vartheta ,t,x\right)$. In this case the function
$u\left(t,x\right)=u\left(t,x,\vartheta \right)$ satisfying the equation
\eqref{3} depends on unknown parameter $\vartheta $ and we can not put
$Y_t=u\left(t,X_t,\vartheta \right)$ because we do not know $\vartheta $.  We
consider the problem of the construction of the couple $(\hat Y_t, \hat Z_t)$,
where $\hat Y_t$ and $\hat Z_t$ are some approximations of
$\left(Y_t,Z_t\right)$. This approximation is done with the help of the
one-step maximum likelihood estimator $\tilde\vartheta $ as $\hat Y_t=
u(t,X_t,\tilde\vartheta )$ and $\hat Z_t=
a\left(t,X_t\right)u'_x(t,X_t,\tilde\vartheta )$. We are interested by a
situation when the error of this approximation is small. One of the
possibilities to have a small error of approximations is in some sense
equivalent to the situation with the small error of estimation of the
parameter $\vartheta $, then from the continuity of the function
$u\left(t,x,\vartheta \right)$ w.r.t. $\vartheta $, we obtain $\hat Y_T \sim
Y_T=\Phi \left(X_T\right)$. The small error of estimation we can have, besides
others, in the situations when $T\rightarrow \infty $ or when $ a \left(\cdot
\right)\rightarrow 0$ (see, e.g., Kutoyants \cite{Kut04} and \cite{Kut94}).
We propose to study this model in the asymptotics of {\it small noise},
i.e. the diffusion coefficient $a\left(t,x\right)^2$ tends to 0. This allows
us to keep the final time $T$ fixed and, what is as well important, this
asymptotics is easier to treat. We show (under regularity conditions) that the
proposed $\hat Y_t$ is close to $Y_t$ for the small values of $\varepsilon $.

We believe that the presented results can be valid (generalized) for
essentially more general, say, nonlinear  models and the conditions
of regularity can be weakened.

\section{Main result}

We consider the following model. The observed diffusion process
$X^T=\left(X_t,0\leq t\leq T\right)$ is
\begin{align}
\label{4}
{\rm d}X_t=S\left(\vartheta ,t,X_t\right){\rm d}t+\varepsilon \sigma
\left(t,X_t\right)\,{\rm d}W_t, \quad X_0, \; 0\leq t\leq T
\end{align}
where $\vartheta \in \Theta =\left(\alpha ,\beta \right)$ and
\begin{align}
\label{5}
&\left|S\left(\vartheta ,t,x_2\right)-S\left(\vartheta ,t,x_1\right)\right|
+\left| \sigma (t,x_2)-\sigma (t,x_2)\right|\leq L\,\left|x_2-x_1\right|,\\
&\left|S\left(\vartheta ,t,x\right)\right|
+\left| \sigma (t,x)\right|\leq L\,\left(1+\left|x\right|\right).
\nonumber
\end{align}
We are given two functions $f\left(t,x,y,z\right)$, $ \Phi \left(x\right)$
and we have to find a couple of stochastic processes $\left(X_t,Z_t,0\leq
t\leq T\right)$ such
that the solution of the equation ({\it backward SDE})
\begin{equation}
\label{6}
{\rm d}Y_t=-f\left(t,X_t,Y_t,Z_t\right)\,{\rm d}t+Z_t\,{\rm d}W_t,\qquad
Y_0,\quad 0\leq t\leq T
\end{equation}
at point $t=T$ satisfies the condition $Y_T=\Phi \left(X_T\right)$.

As the
solution of \eqref{6} is entirely defined by the initial value $Y_0$ and by the
process $Z^T=\left(Z_t,0\leq t\leq T\right)$ we can seek $Y_0, Z^T $, which
provide the equality $ Y_T=\Phi \left(X_T\right)$.

Let us introduce a family of functions
$$
{\cal U}=\left\{\left(u(t,x,\vartheta ),
  t\in \left[0,T\right], x\in {\RR}\right), \vartheta \in \Theta \right\}
$$
 such that for all $\vartheta \in \Theta $ the function $u(t,x,\vartheta ) $
  satisfies the equation
\begin{align*}
\frac{\partial u}{\partial t}+S(\vartheta,t,x)\frac{\partial
  u}{\partial x}
        +\frac{\varepsilon ^2\sigma (t,x)^2}{2} \frac{\partial^2 u}{\partial x^2}
        =-f\left(t,x,u,\varepsilon \sigma (x)\frac{\partial u}{\partial x}\right)
         \end{align*}
and condition $u(T,x,\vartheta )=\Phi \left(x\right)$. If we put
$Y_t=u\left(t,X_t,\vartheta \right)$, then by It\^o's formula we
obtain \eqref{6} with $Z_t=\varepsilon \sigma \left(t,X_t\right)
  u'_x \left(t,X_t,\vartheta \right) $.

We suppose that the true value $\vartheta_0 $ of $\vartheta $ is
unknown. Therefore we can not put $ Y_t=u\left(t,X_t,\vartheta_0 \right)$ and
our goal is to approximate $Y_t$ and $Z_t$. We would like to study this
problem in the situation where the error of approximation can be small.

Of course, the {\it natural approximation} is first to estimate $\vartheta_0 $
and then to substitute it in the function $u\left(\cdot \right)$. The small
error we can have, besides others, in the case when we have a large volume of
observations ($T\rightarrow \infty $) or when the {\it noise} $\varepsilon
\sigma \left(t,X_t\right)$ is small. At the present work we propose an
approximation of $Y_t$ in the case of {\it small noise}, as $\varepsilon
\rightarrow 0$. We suppose to treat {\it large samples case} later.

Remind that the stochastic process  $X_t$ of the equation \eqref{4} under condition
\eqref{5} converges to the deterministic function $x_s=x_s\left(\vartheta
_0\right)$, where $x_s\left(\vartheta \right)$ is  solution of the
ordinary differential equation
\begin{equation}
\label{8}
\frac{{\rm d}x_s}{{\rm d}s}=S\left(\vartheta,s,x_s\right),\qquad x_0, \quad 0\leq s\leq T,
\end{equation}
and this convergence is uniform in $s\in\left[0,T\right]$ (see, e.g.,
\cite{FW} or \cite{Kut94}).
The corresponding PDE with $\varepsilon =0$ is
\begin{align*}
\frac{\partial u^0}{\partial t}+S(\vartheta,t,x)\frac{\partial
  u^0}{\partial x}                 =-f\left(t,x,u^0,0\right),\qquad
u^0\left(T,x,\vartheta \right)=\Phi \left(x\right).
  \end{align*}
and the {\it limit BSDE}
$$
\frac{{\rm d} y_t}{{\rm d}t}=-f\left(t,x_t,y_t,0\right),\qquad y_T=\Phi \left(x_T\right)
$$
we obtain by putting $y_t=u^0\left(t,x_t,\vartheta\right)$.

To estimate $\vartheta $ we can use any {\it good} estimator.
For example, let us denote the likelihood ratio
\begin{align*}
L\left(\vartheta ,X^T\right)=\exp\left\{\int_{0}^{T}\frac{S\left(\vartheta
  ,t,X_t\right)}{\varepsilon ^2\,\sigma \left(t,X_t\right)^2}\,{\rm
  d}X_t-\int_{0}^{T}\frac{S\left(\vartheta ,t,X_t\right)^2}{2\,\varepsilon
  ^2\,\sigma \left(t,X_t\right)^2}\,{\rm d }t\right\}
\end{align*}
and define the maximum likelihood estimator $\hat\vartheta _\varepsilon $ by
the relation
\begin{align*}
L\left(\hat\vartheta_\varepsilon  ,X^T\right)=\sup_{\vartheta \in \Theta
}L\left(\vartheta ,X^T\right).
\end{align*}

Note that we can not use this MLE $\hat\vartheta _\varepsilon $ and to write $\hat
Y_t=u\left(t,X_t,\hat\vartheta
_{\varepsilon} \right)$ because $\hat\vartheta _\varepsilon $
  depends on all observations and  at the moment $t$ the observations $X_s,t< s\leq T$
are not available. If we decide to use just the observations up to instant $t$
and to  define the MLE as follows
\begin{align}
\label{9}
L\left(\hat\vartheta_{t,\varepsilon } ,X^t\right)=\sup_{\vartheta \in \Theta
}L\left(\vartheta ,X^t\right) ,
\end{align}
then we obtain mathematically correct approximation $\hat Y_t=u\left(t,X_t,\hat\vartheta
_{t,\varepsilon} \right)$ and the properties of $\hat Y_t$ are described
in \cite{ZH2013}. Remind that under regularity conditions the
estimator $\hat\vartheta_{t,\varepsilon} $ is consistent, asymptotically
normal
$$
\frac{\hat\vartheta_{t,\varepsilon} -\vartheta_0 }{\varepsilon
}\Longrightarrow {\cal N}\left(0,{\rm I}\left(\vartheta
,x^t\right)^{-1}\right),\qquad {\rm I}\left(\vartheta
,x^t\right)=\int_{0}^{t}\frac{\dot S\left(\vartheta ,s,x_s\right)^2}{\sigma
  \left(s,x_s\right)^2}\,{\rm d}s
$$
and asymptotically efficient (see \cite{Kut94}). Here and in the
sequel dot means derivative w.r.t. $\vartheta $ and  ${\rm
I}\left(\vartheta ,x^t\right)$ is the Fisher information. The
approximation  $\hat Y_t=u\left(t,X_t,\hat\vartheta _{t,\varepsilon}
\right)$ is difficult to realize because to solve the equation
\eqref{9} for all $t\in (0,T]$ is computationally a quite complicate
problem.

We need some regularity conditions. Let us denote ${\cal P}$ a class
of functions of $x$ having polynomial majorants. For example a
function $g\left(t,x,\vartheta ,\varepsilon \right)\in {\cal P}$
means that there exist constants $C>0$ and $p>0$ which do not depend
on $t\in \left[0,T\right],\vartheta\in   \Theta
,\varepsilon\in\left[0,1\right]$ such that
\begin{equation}
\label{7}
\left|g\left(t,x,\vartheta ,\varepsilon \right)\right|\leq C\left(1+\left|x\right|^p\right) .
\end{equation}

 We suppose that  the functions $S\left(\vartheta ,t,x\right)$ and  $u(t,x,\vartheta )   $  have two
continuous derivatives w.r.t. $\vartheta  $ and the following derivatives
belong to ${\cal P} $
 $$
\dot S\left(\vartheta
,t,x\right),\ddot S\left(\vartheta ,t,x\right),\dot S'_x\left(\vartheta
,t,x\right), \sigma'_x \left(t,x\right), \dot  u(t,x,\vartheta ),\ddot
u(t,x,\vartheta ),\dot  u'_x(t,x,\vartheta )  .
$$
The function $\sigma \left(t,x\right)^2\geq \kappa >0$ and we have the uniform
in $\vartheta $
convergence (as $\varepsilon \rightarrow 0$)
\begin{equation}
\label{q}
u\left(t,x,\vartheta \right)\rightarrow u^0\left(t,x,\vartheta \right),\qquad u'_x\left(t,x,\vartheta \right)\rightarrow\left( u^0\right)'_x\left(t,x,\vartheta \right).
\end{equation}

We propose the following solution. Fix some (small) $\delta >0$ and introduce the {\it
  minimum distance estimator } (MDE) $\vartheta ^*_{\delta ,\varepsilon }$ by the relation
\begin{align*}
\left\|X-x\left(\vartheta ^*_{\delta ,\varepsilon
}\right)\right\|^2=\inf_{\vartheta \in \Theta }\left\|X-x\left(\vartheta
\right)\right\|^2=\inf_{\vartheta \in \Theta }\int_{0}^{\delta
}\left[X_t-x_t\left(\vartheta \right)\right] ^2\,{\rm d} t.
\end{align*}

This estimator  is consistent and asymptotically
normal
$$
\varepsilon ^{-1}\left( \vartheta ^*_{\delta ,\varepsilon }-\vartheta _0
\right)\Longrightarrow {\cal N} \left(0,D_\delta \left(\vartheta
_0\right)^2\right)
$$
where $D_\delta \left(\vartheta _0\right)^2>0 $ (see Theorem 7.5 \cite{Kut94}.
The required regularity conditions are  : the function $S\left(\vartheta
,t,x\right)$ has two continuous derivatives w.r.t. $\vartheta $ having
polynomial majorants (see \eqref{7}) and the following identifiability condition is fulfilled:
for any $\nu >0$
$$
\inf_{\left|\vartheta -\vartheta _0\right|>\nu }\left\|x\left(\vartheta
\right)-x\left(\vartheta_0 \right)\right\| >0.
$$
Let us introduce the {\it one-step MLE}
\begin{equation}
\label{10}
\tilde\vartheta _{t,\varepsilon }=\vartheta ^*_{\delta ,\varepsilon
}+\frac{\Delta _t\left(\vartheta ^*_{\delta ,\varepsilon },X_\delta
  ^t\right)+\Delta _\delta \left(\vartheta ^*_{\delta ,\varepsilon
  },X^\delta\right)}{{\rm
    I}\left(\vartheta ^*_{\delta ,\varepsilon } ,x^t\left(\vartheta ^*_{\delta
    ,\varepsilon }\right)\right) } ,
\end{equation}
where
\begin{align*}
&\Delta _t\left(\vartheta,X_\delta ^t\right)=\int_{\delta }^{t}\frac{\dot
  S\left(\vartheta ,s,X_s\right)}{ \sigma
  \left(s,X_s\right)^2}\,\left[{\rm d}X_s-S\left(\vartheta ,s,X_s\right)\,{\rm
    d}s\right],\quad t\in [\delta ,T],\\
&\Delta _\delta \left(\vartheta,X^\delta\right)=A\left(\vartheta
,\delta ,X_\delta \right)-\int_{0}^{\delta }A'_s\left(\vartheta
,s ,X_s \right)\,{\rm d}s\\
& \quad  -\frac{\varepsilon^2 }{2}\int_{0}^{\delta } B'_x\left(\vartheta ,s
,X_s\right)\sigma \left(s,X_s\right)^2{\rm d}s-\int_{0 }^{\delta }\frac{\dot
  S\left(\vartheta ,s,X_s\right)S\left(\vartheta ,s,X_s\right)}{  \sigma
  \left(s,X_s\right)^2}{\rm d}s ,\\
&B\left(\vartheta ,s,x\right)=\frac{\dot
  S\left(\vartheta ,s,x\right)}{ \sigma
  \left(s,x\right)^2},\qquad A\left(\vartheta
,s,x\right)=\int_{x_0}^{x}B\left(\vartheta ,s,z\right)\,{\rm d}z, \\
&{\rm I}\left(\vartheta
  ,x^t\left(\vartheta\right)\right)=\int_{0}^{t}\frac{\dot S\left(\vartheta
    ,s,x_s\left(\vartheta \right)\right)^2}{\sigma \left(s,x_s\left(\vartheta
    \right)\right)^2}\,{\rm d}s.
\end{align*}

The approximation of the solution of BSDE  is given in  the following theorem.
\begin{theorem}
\label{T1} Let the conditions of regularity be fulfilled then the processes
$$
\hat Y_t=u\left(t,X_t,\tilde\vartheta _{t,\varepsilon }\right),\qquad
\hat Z_t=\varepsilon \sigma \left(t,X_t\right)  u'_x\left(t,X_t,\tilde\vartheta
_{t,\varepsilon }\right)
$$
for the values $t\in \left[\delta ,T\right]$ have the representation
\begin{align}
\label{11}
\hat Y_t&=Y_t+\varepsilon \dot u\left(t,X_t,\vartheta _0\right)\,\xi _t\left(\vartheta
_0\right)+o\left(\varepsilon \right) ,\\
\hat Z_t&=Z_t+\varepsilon ^2\sigma \left(t,X_t\right)\dot u'_x\left(t,X_t,\vartheta _0\right)\,\xi
_t\left(\vartheta _0\right)+o\left(\varepsilon^2 \right)  ,
\label{12}
\end{align}
where
$$
\xi _t\left(\vartheta_0\right)={\rm
    I}\left(\vartheta_{0} ,x^t\left(\vartheta _{0}\right)\right)^{-1} \int_{0}^{t}\frac{\dot
  S\left(\vartheta_0 ,s,x_s\left(\vartheta_0\right)\right)}{\sigma
  \left(s,x_s\left(\vartheta_0\right)\right)}\,{\rm d}W_s.
$$
\end{theorem}
{\bf Proof.}
Suppose that we already proved that
\begin{equation}
\label{13}
\varepsilon ^{-1}\left(\tilde \vartheta _{t,\varepsilon }-\vartheta
_0\right)= \xi _t\left(\vartheta
_0\right)+o\left(1\right),
\end{equation}
then the representations \eqref{11}, \eqref{12}  we obtain by
Taylor's formula
\begin{align*}
\hat Y_t&= u\left(t,X_t, \vartheta _{0 }\right)+ \left(\tilde \vartheta
_{t,\varepsilon }-\vartheta _0\right) \dot u\left(t,X_t, \vartheta _{0
}\right)+ o\left(\varepsilon \right)\\
&=Y_t+\varepsilon \dot u\left(t,X_t, \vartheta _{0
}\right)\xi _t\left(\vartheta _0\right)+ o\left(\varepsilon \right) ,
\end{align*}
and
\begin{align*}
\hat Z_t&=\varepsilon \sigma
\left(t,X_t\right) u'_x\left(t,X_t, \vartheta _{0 }\right)+ \left(\tilde
\vartheta _{t,\varepsilon }-\vartheta _0\right)\varepsilon \sigma
\left(t,X_t\right) \dot u'_x\left(t,X_t, \vartheta
_{0 }\right)+ o\left(\varepsilon \right)\\
&=Z_t+\varepsilon ^2\sigma
\left(t,X_t\right) \dot u'_x\left(t,X_t, \vartheta
_{0 }\right)\xi _t\left(\vartheta _0\right)+ o\left(\varepsilon^2 \right).
\end{align*}
Let us verify \eqref{13}. Remind that for any $p>0$
\begin{align}
&\sup_{0\leq t\leq T}\left|X_t-x_t\left(\vartheta _0\right)\right|\leq
 C\varepsilon\,\sup_{0\leq t\leq T}\left|W_t\right|,\nonumber\\
&\Ex_{\vartheta _0}\left|X_t-x_t\left(\vartheta _0\right)\right|^p\leq
 C\,\varepsilon ^p
\label{14}\\
&\Ex_{\vartheta _0}\left|\vartheta ^*_{\delta  ,\varepsilon }-\vartheta
_0\right|^p\leq C\,\varepsilon ^{p} \nonumber
\end{align}
(see, e.g.; Lemma 1.13, and Theorem 7.5 in  \cite{Kut94}).
Below we denoted $x_s=x_s\left(\vartheta _0\right)$, use the convergence
\eqref{14}, consistency of the MDE $\vartheta ^*_{\delta ,\varepsilon
} $ and smoothness of $S\left(\vartheta ,s,x\right)$ and $\sigma
\left(s,x\right)$
\begin{align*}
&\Delta _t\left(\vartheta ^*_{\delta ,\varepsilon
},X_\delta ^t\right)=\int_{\delta }^{t}\frac{\dot
  S\left(\vartheta ^*_{\delta ,\varepsilon
} ,s,X_s\right)}{ \sigma
  \left(s,X_s\right)^2}\,\left[{\rm d}X_s-S\left(\vartheta ^*_{\delta ,\varepsilon
} ,s,X_s\right)\,{\rm     d}s\right]\\
&\qquad =\varepsilon\int_{\delta }^{t}\frac{\dot
  S\left(\vartheta ^*_{\delta ,\varepsilon
} ,s,X_s\right)}{\sigma
  \left(s,X_s\right)}\,{\rm d}W_s\\
&\qquad \qquad +\int_{\delta }^{t}\frac{\dot
  S\left(\vartheta ^*_{\delta ,\varepsilon
} ,s,X_s\right)}{ \sigma
  \left(s,X_s\right)^2}\,\left[S\left(\vartheta _{0
} ,s,X_s\right)-S\left(\vartheta ^*_{\delta ,\varepsilon
} ,s,X_s\right)\right]\,{\rm     d}s\\
&\qquad =\varepsilon\int_{\delta }^{t}\frac{\dot
  S\left(\vartheta_{0
} ,s,x_s\right)}{\sigma
  \left(s,x_s\right)}\,{\rm d}W_s-{\left(\vartheta ^*_{\delta ,\varepsilon
}-\vartheta _0\right)}{ }\int_{\delta }^{t}\frac{\dot
  S\left(\vartheta_{0
} ,s,x_s\right)^2}{ \sigma
  \left(s,x_s\right)^2}\,\,{\rm     d}s+o\left(\varepsilon\right).
\end{align*}
Further, using the same arguments we write
\begin{align*}
 &\Delta _\delta \left(\vartheta ^*_{\delta ,\varepsilon
},X^\delta\right)=A\left(\vartheta ^*_{\delta ,\varepsilon } ,\delta
,X_\delta \right)-\int_{0}^{\delta }A'_s\left(\vartheta ^*_{\delta
  ,\varepsilon } ,s ,X_s \right)\,{\rm d}s\\
& \quad
-\frac{\varepsilon^2 }{2}\int_{0}^{\delta } B'_x\left(\vartheta ^*_{\delta
  ,\varepsilon } ,s ,X_s\right)\sigma \left(s,X_s\right)^2{\rm d}s-\int_{0 }^{\delta }\frac{\dot
  S\left(\vartheta ^*_{\delta ,\varepsilon  },s,X_s\right)S\left(\vartheta
  ^*_{\delta ,\varepsilon } ,s,X_s\right)}{  \sigma
  \left(s,X_s\right)^2}{\rm d}s \\
&\qquad \quad =A\left(\vartheta_{0 } ,\delta
,X_\delta \right)-\int_{0}^{\delta }A'_s\left(\vartheta _{0} ,s ,X_s \right)\,{\rm d}s\\
&-\frac{\varepsilon^2 }{2}\int_{0}^{\delta } B'_x\left(\vartheta_{0 }  ,s
,X_s\right)\sigma \left(s,X_s\right)^2{\rm d}s-\int_{0 }^{\delta }\frac{\dot
  S\left(\vartheta_{0 } ,s,X_s\right)S\left(\vartheta_{0 } ,s,X_s\right)}{
  \sigma   \left(s,X_s\right)^2}{\rm d}s \\
&\quad-\left(\vartheta ^*_{\delta ,\varepsilon }-\vartheta_{0} \right) \int_{0 }^{\delta }\frac{\dot
  S\left(\vartheta_{0 } ,s,X_s\right)^2}{
  \sigma   \left(s,X_s\right)^2}{\rm d}s+ \left(\vartheta ^*_{\delta
  ,\varepsilon }-\vartheta ^*_{0} \right) \left[\dot A\left(\vartheta _0 ,\delta
,X_\delta \right)\right.\\
&\quad -\int_{0}^{\delta } \dot A'_s\left(\vartheta_0 ,s ,X_s \right)\,{\rm
    d}s-\frac{\varepsilon^2 }{2}\int_{0}^{\delta }\dot  B'_x\left(\vartheta_0
  ,s ,X_s\right)\sigma \left(s,X_s\right)^2{\rm d}s\\
&\quad\left. -\int_{0 }^{\delta }\frac{\ddot
  S\left(\vartheta_{0 } ,s,X_s\right)S\left(\vartheta_{0 } ,s,X_s\right)}{
  \sigma   \left(s,X_s\right)^2}{\rm d}s\right]+o\left(\varepsilon \right).
\end{align*}
Note that by It\^o's formula
\begin{align*}
&A\left(\vartheta_{0 } ,\delta
,X_\delta \right)-\int_{0}^{\delta }A'_s\left(\vartheta _{0} ,s ,X_s
  \right)\,{\rm d}s-\frac{\varepsilon^2 }{2}\int_{0}^{\delta }
  B'_x\left(\vartheta_{0 }  ,s
,X_s\right)\sigma \left(s,X_s\right)^2{\rm d}s\\
&-\int_{0 }^{\delta }\frac{\dot
  S\left(\vartheta_{0 } ,s,X_s\right)S\left(\vartheta_{0 } ,s,X_s\right)}{
  \sigma   \left(s,X_s\right)^2}{\rm d}s =\varepsilon \int_{0}^{\delta }\frac{\dot
  S\left(\vartheta_{0 } ,s,X_s\right)}{
  \sigma   \left(s,X_s\right)}{\rm d}W_s
\end{align*}
and
\begin{align*}
&\dot A\left(\vartheta_{0 } ,\delta
,X_\delta \right)-\int_{0}^{\delta }\dot A'_s\left(\vartheta _{0} ,s ,X_s
  \right)\,{\rm d}s-\frac{\varepsilon^2 }{2}\int_{0}^{\delta }\dot
  B'_x\left(\vartheta_{0 }  ,s
,X_s\right)\sigma \left(s,X_s\right)^2{\rm d}s\\
&\qquad -\int_{0 }^{\delta }\frac{\ddot
  S\left(\vartheta_{0 } ,s,X_s\right)S\left(\vartheta_{0 } ,s,X_s\right)}{
  \sigma   \left(s,X_s\right)^2}{\rm d}s =\varepsilon \int_{0}^{\delta }\frac{\ddot
  S\left(\vartheta_{0 } ,s,X_s\right)}{
  \sigma   \left(s,X_s\right)}{\rm d}W_s.
\end{align*}
Therefore
\begin{align*}
&\Delta _t \left(\vartheta ^*_{\delta ,\varepsilon
},X^t_\delta\right)+\Delta _\delta \left(\vartheta ^*_{\delta ,\varepsilon
},X^\delta\right)=\varepsilon \int_{0}^{t }\frac{\dot
  S\left(\vartheta_{0 } ,s,X_s\right)}{
  \sigma   \left(s,X_s\right)}{\rm d}W_s\\
&\qquad \qquad - \left(\vartheta ^*_{\delta ,\varepsilon }-\vartheta_{0} \right)  \int_{0}^{t }\frac{\dot
  S\left(\vartheta_{0 } ,s,X_s\right)^2}{
  \sigma   \left(s,X_s\right)^2}{\rm d}s+o\left(\varepsilon \right)\\
&\qquad =\varepsilon \int_{0}^{t }\frac{\dot
  S\left(\vartheta_{0 } ,s,x_s\right)}{
  \sigma   \left(s,x_s\right)}{\rm d}W_s- \left(\vartheta ^*_{\delta
    ,\varepsilon }-\vartheta_{0} \right)  \int_{0}^{t }\frac{\dot
  S\left(\vartheta_{0 } ,s,x_s\right)^2}{
  \sigma   \left(s,x_s\right)^2}{\rm d}s+o\left(\varepsilon \right).
\end{align*}
For the one-step MLE this allows us to write
\begin{align*}
&\varepsilon ^{-1}\left(\tilde\vartheta _{t,\varepsilon }-\vartheta _0\right)
=\frac{\vartheta^* _{\delta ,\varepsilon }-\vartheta
_0}{\varepsilon } +{\rm I}\left(\vartheta^* _{\delta ,\varepsilon
},x^t\left(\vartheta^* _{\delta ,\varepsilon }
\right)\right)^{-1}\left[\int_{0}^{t }\frac{\dot
  S\left(\vartheta_{0 } ,s,x_s\right)}{
  \sigma   \left(s,x_s\right)}{\rm d}W_s\right.\\
&\qquad \quad \left.\qquad - \frac{\left(\vartheta ^*_{\delta
    ,\varepsilon }-\vartheta_{0} \right)}{\varepsilon }  \int_{0}^{t }\frac{\dot
  S\left(\vartheta_{0 } ,s,x_s\right)^2}{
  \sigma   \left(s,x_s\right)^2}{\rm d}s\right]+o\left(1\right)
  \\
&\qquad = {\rm I}\left(\vartheta _0,x^t\left(\vartheta _0 \right)\right)^{-1}\int_{0}^{t }\frac{\dot
  S\left(\vartheta_{0 } ,s,x_s\right)}{
  \sigma   \left(s,x_s\right)}{\rm d}W_s+o\left(1\right)=\xi _t\left(\vartheta _0\right)+o\left(1\right).
\end{align*}
This proves Theorem \ref{T1}. We used just the continuity of the derivatives
$\dot u$

Let us show that the proposed approximation is asymptotically
efficient. This means, that the means-quare errors
$$
\Ex_\vartheta \left|Y_t-\hat Y_t\right|^2,\qquad \Ex_\vartheta \left|Z_t-\hat Z_t\right|^2,
$$
of estimation $Y_t$ and $Z_t$  can not be improved. This will be
done in two steps. First we establish a low bound on the risks of
all estimators and then show that the proposed estimators attaint
this bound.

\begin{theorem}
\label{T2} For all estimators $\bar Y_t$ and $\bar Z_t$ and all $t\in
\left[\delta ,T\right]$ we have the relations
\begin{align}
\label{15}
&\Liminf_{\nu \rightarrow 0}\Liminf_{\varepsilon \rightarrow 0} \sup_{\left|\vartheta
  -\vartheta _0\right|\leq \nu } \varepsilon ^{-2}\Ex_\vartheta \left| \bar
Y_t-Y_t\right|^2\geq \frac{\dot u^0\left(t,x_t\left(\vartheta
  _0\right),\vartheta _0\right)^2}{{\rm I}\left(\vartheta
  _0,x^t\left(\vartheta _0 \right)\right)} ,\\
&\Liminf_{\nu \rightarrow 0}\Liminf_{\varepsilon \rightarrow 0} \sup_{\left|\vartheta
  -\vartheta _0\right|\leq \nu } \varepsilon ^{-4}\Ex_\vartheta \left| \bar
Z_t-Z_t\right|^2\geq \frac{\left(\dot u^{ 0}\right)'_x\left(t,x_t\left(\vartheta
  _0\right),\vartheta _0\right)^2\sigma \left(t,x_t\left(\vartheta
  _0\right)\right)^2}{{\rm I}\left(\vartheta   _0,x^t\left(\vartheta _0
  \right)\right)}
\label{16}
\end{align}
\end{theorem}
{\bf Proof.} We follow the usual proof of the van Trees and minimax
bounds. Let us fix $\nu >0$ and introduce a probability density $p\left(\theta
\right),\theta \in \left[\vartheta _0-\nu ,\vartheta _0+\nu \right]$ such that
$p\left(\vartheta _0-\nu \right)=0$, $p\left(\vartheta _0+\nu \right)=0$  and
$$
{\rm I}_p=\int_{\vartheta _0-\nu}^{\vartheta _0+\nu}\frac{\dot p\left(\theta
  \right)^2}{p\left(\theta \right)}\,{\rm d}\theta  <\infty .
$$
Denote $L_{\theta _0}\left(\theta ,X^t\right)=L\left(\theta_0
,X^t\right)^{-1}L\left(\theta ,X^t\right)$.  Integrating by parts we obtain
\begin{align*}
\int_{\vartheta _0-\nu}^{\vartheta _0+\nu}& \left[\bar Y_t- u\left(t,X_t,\vartheta
\right)\right]\frac{\partial }{\partial \theta }\left[L_{\theta _0}\left(\theta , X^t\right)p\left(\theta
\right) \right]{\rm d} \theta\\
&\qquad=\left.L_{\theta _0}\left(\theta ,X^t\right)p\left(\theta
\right)\left[\bar Y_t-u\left(t,X_t,\vartheta
\right)\right]\right|_{\vartheta _0-\nu}^{\vartheta _0+\nu}\\
&\qquad\qquad \quad  +\int_{\vartheta
  _0-\nu}^{\vartheta _0+\nu}\dot u\left(t,X_t,\vartheta\right)L_{\theta _0}\left(\theta , X^t\right)p\left(\theta\right){\rm d}\theta \\
&\qquad =\int_{\vartheta
  _0-\nu}^{\vartheta _0+\nu}\dot u\left(t,X_t,\vartheta\right)L_{\theta _0}\left(\theta ,X^t\right)p\left(\theta
\right){\rm d}\theta .
\end{align*}
We have
\begin{align*}
&\Ex_{\vartheta _0} \int_{\vartheta _0-\nu}^{\vartheta _0+\nu} \left[\bar Y_t- u\left(t,X_t,\vartheta
\right)\right]\frac{\partial }{\partial \theta }\left[L_{\theta _0}\left(\theta , X^t\right)p\left(\theta
\right) \right]{\rm d} \theta\\
&\qquad =\Ex_{\vartheta _0} \int_{\vartheta
  _0-\nu}^{\vartheta _0+\nu}\dot u\left(t,X_t,\vartheta\right)L_{\theta _0}\left(\theta ,X^t\right)p\left(\theta
\right){\rm d}\theta \\
&\qquad =\int_{\vartheta
  _0-\nu}^{\vartheta _0+\nu}\Ex_{\vartheta}\dot u\left(t,X_t,\vartheta\right)p\left(\theta
\right){\rm d}\theta=\int_{\vartheta
  _0-\nu}^{\vartheta _0+\nu}\dot u\left(t,x_t,\vartheta\right)p\left(\theta
\right){\rm d}\theta+O\left(\varepsilon \right) .
\end{align*}
We need the equality
\begin{align*}
&\Ex_{\vartheta _0} \int_{\vartheta _0-\nu}^{\vartheta
    _0+\nu}\left(\frac{\partial }{\partial \theta }\ln \left[L_{\theta
      _0}\left(\theta , X^t\right)p\left(\theta \right)
    \right]\right)^2L_{\theta _0}\left(\theta , X^t\right)p\left(\theta
  \right) {\rm d} \theta\\
&\quad = \int_{\vartheta _0-\nu}^{\vartheta
    _0+\nu}\Ex_{\vartheta } \left(\int_{0}^{t}\frac{\dot
    S\left(\vartheta ,s,X_s\right)}{\varepsilon \;\sigma \left(s,X_s\right)}{\rm
    d}W_s+\frac{\dot p\left(\theta \right)}{p\left(\theta \right)} \right)^2
  p\left(\theta \right) {\rm d} \theta\\
&\quad = \int_{\vartheta _0-\nu}^{\vartheta _0+\nu}\Ex_{\theta
  } \int_{0}^{t}\frac{\dot S\left(\theta ,s,X_s\right)^2}{\varepsilon
  ^2\sigma
    \left(s,X_s\right)^2}{\rm d}s\; p\left(\theta \right) \;{\rm d}\theta +  {\rm     I}_p\\
&\quad =\frac{1}{ \varepsilon  ^2}\int_{\vartheta _0-\nu}^{\vartheta _0+\nu} {\rm     I}\left(\theta,x^t\left(\theta \right)
  \right)\; p\left(\theta \right) \;{\rm d}\theta +   {\rm     I}_p+o\left(\frac{1}{\varepsilon}\right),
\end{align*}
where we used  the chain
rule  ($\Ex_{\vartheta _0}L_{\theta _0}\left(\theta ,
X^t\right)=\Ex_{\vartheta} $).
Below we apply Cauchy-Shwartz inequality
\begin{align*}
&\left(\Ex_{\vartheta _0} \int_{\vartheta _0-\nu}^{\vartheta _0+\nu} \left[\bar Y_t- u\left(t,X_t,\vartheta
\right)\right]\frac{\partial }{\partial \theta }\left[L_{\theta
      _0}\left(\theta , X^t\right)p\left(\theta \right)\right]{\rm d}
  \theta\right)^2\\
&\quad =
\left(\Ex_{\vartheta _0} \int_{\vartheta _0-\nu}^{\vartheta _0+\nu} \left[\bar Y_t- u\left(t,X_t,\vartheta
\right)\right]\right.\\
&\qquad \qquad\qquad \qquad\left. \frac{\partial }{\partial \theta }\ln
\left[L_{\theta _0}\left(\theta , X^t\right)p\left(\theta \right)
  \right]L_{\theta _0}\left(\theta , X^t\right)p\left(\theta
\right) {\rm d} \theta\right)^2 \\
&\quad \leq  \Ex_{\vartheta _0} \int_{\vartheta _0-\nu}^{\vartheta _0+\nu}
\left[\bar Y_t- u\left(t,X_t,\vartheta
\right)\right]^2L_{\theta _0}\left(\theta , X^t\right)p\left(\theta
\right) {\rm d} \theta\\
&\qquad \qquad \Ex_{\vartheta _0} \int_{\vartheta _0-\nu}^{\vartheta
  _0+\nu}\left(\frac{\partial }{\partial \theta }\ln \left[L_{\theta
    _0}\left(\theta , X^t\right)p\left(\theta \right)
  \right]\right)^2L_{\theta _0}\left(\theta , X^t\right)p\left(\theta
\right) {\rm d} \theta\\
&= \int_{\vartheta _0-\nu}^{\vartheta _0+\nu} \Ex_{\vartheta }\left[\bar Y_t- u\left(t,X_t,\vartheta
\right)\right]^2p\left(\theta\right) {\rm d} \theta\\
&\qquad \qquad
\left[\frac{1}{ \varepsilon  ^2}\int_{\vartheta _0-\nu}^{\vartheta _0+\nu} {\rm     I}\left(\theta,x^t\left(\theta \right)
  \right)\; p\left(\theta \right) \;{\rm d}\theta +   {\rm     I}_p+o\left(\frac{1}{\varepsilon }\right)\right].
\end{align*}
Therefore we obtained the van Trees inequality
\begin{align*}
&\int_{\vartheta _0-\nu}^{\vartheta _0+\nu} \Ex_{\vartheta }\left[\bar Y_t- u\left(t,X_t,\vartheta
\right)\right]^2p\left(\theta\right) {\rm d} \theta\\
&\qquad \geq \frac{\left(\int_{\vartheta
  _0-\nu}^{\vartheta _0+\nu}\dot u\left(t,x_t,\vartheta\right)p\left(\theta
\right){\rm d}\theta+O\left(\varepsilon \right) \right)^2 }{ \frac{1}{\varepsilon  ^2}\int_{\vartheta _0-\nu}^{\vartheta _0+\nu} {\rm     I}\left(\theta,x^t\left(\theta \right)
  \right)\; p\left(\theta \right) \;{\rm d}\theta +   {\rm     I}_p+o\left(\frac{1}{\varepsilon }\right)\  }
\end{align*}
Further
\begin{align*}
&\Liminf_{\varepsilon \rightarrow 0}\sup_{\left|\vartheta -\vartheta
    _0\right|\leq \nu }\varepsilon ^{-2}\Ex_\vartheta \left| \bar
  Y_t-Y_t\right|^2=\Liminf_{\varepsilon \rightarrow 0}\sup_{\left|\vartheta
    -\vartheta _0\right|\leq \nu }\varepsilon ^{-2}\Ex_\vartheta \left| \bar Y_t-
  u\left(t,X_t,\vartheta \right)\right|^2\\ &\quad \geq \Liminf_{\varepsilon
    \rightarrow 0}\varepsilon ^{-2}\int_{\vartheta _0-\nu}^{\vartheta _0+\nu}
  \Ex_{\vartheta }\left[\bar Y_t- u\left(t,X_t,\vartheta
    \right)\right]^2p\left(\theta\right) {\rm d} \theta\\ &\quad \geq
  \frac{\left(\int_{\vartheta _0-\nu}^{\vartheta _0+\nu}\dot
    u^0\left(t,x_t,\vartheta\right)p\left(\theta \right){\rm d}\theta \right)^2
  }{ \int_{\vartheta _0-\nu}^{\vartheta _0+\nu} {\rm
      I}\left(\theta,x^t\left(\theta \right) \right)\; p\left(\theta \right)
    \;{\rm d}\theta }\longrightarrow \frac{\dot u^0\left(t,x_t\left(\vartheta_0
    \right),\vartheta_0\right)^2 }{  {\rm
      I}\left(\vartheta_0,x^t\left(\vartheta_0 \right) \right)}.
\end{align*}
The last limit corresponds to $\nu \rightarrow 0$ and used the continuity of
the underlying functions. Therefore the inequality \eqref{15} is proved. The
proof of \eqref{16} is quite similar.

We call an approximation $Y_t^\star$ asymptotically efficient if for
all $\vartheta _0\in \Theta $ we have the equality
\begin{equation}
\label{17}
\lim_{\nu \rightarrow 0}\lim_{\varepsilon \rightarrow 0} \sup_{\left|\vartheta
  -\vartheta _0\right|\leq \nu } \varepsilon ^{-2}\Ex_\vartheta \left|
Y_t^\star-Y_t\right|^2= \frac{\dot u^0\left(t,x_t\left(\vartheta
  _0\right),\vartheta _0\right)^2}{{\rm I}\left(\vartheta
  _0,x^t\left(\vartheta _0 \right)\right)}
\end{equation}
and of course the similar definition is valid in the case of the bound \eqref{16}.

\begin{theorem}
\label{T3} The approximations $$\hat Y_t=u\left(t,X_t,\tilde \vartheta
_{t,\varepsilon }\right)\quad {\rm  and } \quad \hat Z_t=\varepsilon \sigma
\left(t,X_t\right)u'_x\left(t,X_t,\tilde \vartheta _{t,\varepsilon }\right)$$
are asymptotically efficient, i.e.,
\begin{align*}
\lim_{\nu \rightarrow 0}\lim_{\varepsilon \rightarrow 0} \sup_{\left|\vartheta
  -\vartheta _0\right|\leq \nu } \varepsilon ^{-2}\Ex_\vartheta \left|
\hat Y_t-Y_t\right|^2&= \frac{\dot u^0\left(t,x_t\left(\vartheta
  _0\right),\vartheta _0\right)^2}{{\rm I}\left(\vartheta
  _0,x^t\left(\vartheta _0 \right)\right)} ,\\
\lim_{\nu \rightarrow 0}\lim_{\varepsilon \rightarrow 0} \sup_{\left|\vartheta
  -\vartheta _0\right|\leq \nu } \varepsilon ^{-4}\Ex_\vartheta \left|
\hat Z_t-Z_t\right|^2&= \frac{\sigma \left(t,x_t\left(\vartheta
  _0\right)\right)^2\left(\dot u^0\right)'_x\left(t,x_t\left(\vartheta
  _0\right),\vartheta _0\right)^2}{{\rm I}\left(\vartheta
  _0,x^t\left(\vartheta _0 \right)\right)}
\end{align*}
\end{theorem}
{\bf Proof.} The proof of these equalities follows from
\eqref{11}-\eqref{12} because using standard arguments we can show that  the
convergence $o\left(1\right)$ in  these representations is uniform in
$\vartheta $ and the moments converge too.

\section{Discussion}

{\it Uniform approximation.} The representations \eqref{11},
\eqref{12} are valid for each $t\in \left[\delta ,T\right]$. It is
possible to show that these equalities are true uniformly in $t$.
More precisely, let us put $\nu=\varepsilon ^\kappa ,\kappa  >0$.
Then for sufficiently small $\kappa $ we have the convergence
\begin{align*}
\Pb_{\vartheta_0 }\left\{\sup_{\delta \leq t\leq T}\left|\hat Y_t-Y_t\right|>\nu \right\}\longrightarrow 0.
\end{align*}

Indeed, as the derivatives $\dot u\left(t,x,\vartheta \right)$ and
$\dot u'_x\left(t,x,\vartheta \right)$ have polynomial majorants, we
can write
\begin{align*}
&\Pb_{\vartheta_0 }\left\{\sup_{\delta \leq t\leq T}\left|\hat
Y_t-Y_t\right|>\nu \right\}=\Pb_{\vartheta_0 }\left\{\sup_{\delta \leq t\leq
  T}\left|\dot u\left(t,X_t,\bar \vartheta \right)\right|\;\left|
\tilde\vartheta _{t,\varepsilon }-\vartheta _0\right|>\nu \right\}  \\
&\quad \leq \Pb_{\vartheta_0 }\left\{\sup_{\delta \leq t\leq
  T}\left|\dot u\left(t,X_t,\bar \vartheta \right)\right|>\nu^{-\frac{1}{2}} \right\}+\Pb_{\vartheta_0 }\left\{\sup_{\delta \leq t\leq
  T}\left|
\tilde\vartheta _{t,\varepsilon }-\vartheta _0\right|>\nu^{\frac{3}{2}} \right\} \\
&\quad \leq \Pb_{\vartheta_0 }\left\{C\sup_{\delta \leq t\leq
  T}\left|X_t\right|^p>\nu^{-\frac{1}{2}} \right\}+\Pb_{\vartheta_0 }\left\{\sup_{\delta \leq t\leq
  T}\left|
\tilde\vartheta _{t,\varepsilon }-\vartheta _0\right|>\nu^{\frac{3}{2}} \right\}.
\end{align*}
The estimates \eqref{14} allow us to prove    the convergence
$$
\Pb_{\vartheta_0 }\left\{\sup_{\delta \leq t\leq
  T}\left|X_t\right|^p>\nu^{-\frac{1}{2}} \right\}\rightarrow 0.
$$
 The verification
$$
\Pb_{\vartheta_0 }\left\{\sup_{\delta \leq t\leq
  T}\left|
\tilde\vartheta _{t,\varepsilon }-\vartheta _0\right|>\nu^{\frac{3}{2}} \right\}
$$
with $3\kappa <2$ is more complicate, but direct, because we have the uniform
convergence
\begin{align*}
&\sup_{\delta \leq t\leq
  T}\left|{\rm I}\left(\vartheta^* _{\delta ,\varepsilon
},x^t\left(\vartheta^* _{\delta ,\varepsilon }\right)\right)-{\rm
  I}\left(\vartheta _{0 },x^t\left(\vartheta _{0 }\right)\right)\right|
\rightarrow 0 ,\\
&\sup_{\delta \leq t\leq  T}\left|{\rm I}\left(\vartheta^* _{\delta
  ,\varepsilon},X^t\right)-{\rm I}\left(\vartheta _{0},x^t\left(\vartheta _{0
}\right)\right) \right| \rightarrow 0 .
\end{align*}
Further ${\rm I}\left(\vartheta^* _{\delta ,\varepsilon
},x^t\left(\vartheta^* _{\delta ,\varepsilon }\right)\right)\geq {\rm
  I}\left(\vartheta^* _{\delta ,\varepsilon
},x^\delta \left(\vartheta^* _{\delta ,\varepsilon }\right)\right)\geq
\inf_\vartheta {\rm I}\left(\vartheta,x^\delta
 \left(\vartheta\right)\right)>0 $ and
\begin{align*}
&\Pb_{\vartheta_0 }\left\{\sup_{\delta \leq t\leq
  T}\left| \int_{\delta }^{t} \left[\frac{\dot S\left(\vartheta^* _{\delta ,\varepsilon
},s,X_s\right)}{\sigma \left(s,X_s\right)} - \frac{\dot S\left(\vartheta _{0
},s,x_s\right)}{\sigma \left(s,x_s\right)} \right]{\rm d}W_t
  \right|>\nu\right\}\\
&\qquad \leq C\nu ^{-2}\Ex_{\vartheta_0 }\int_{\delta }^{T} \left[\frac{\dot
      S\left(\vartheta^* _{\delta ,\varepsilon
},s,X_s\right)}{\sigma \left(s,X_s\right)} - \frac{\dot S\left(\vartheta _{0
},s,x_s\right)}{\sigma \left(s,x_s\right)} \right]^2{\rm d}s\leq C\nu ^{-2}\varepsilon ^2.
\end{align*}
The proof of the last estimate is direct.

\bigskip

{\it Case $\delta \rightarrow 0$.} The representations \eqref{11}, \eqref{12} are
valid for each $t\in \left[\delta ,T\right]$ with fixed $\delta >0$. It is
possible to show that $\hat Y_t\rightarrow Y_t$ and $\varepsilon ^{-1}\hat
Z_t\rightarrow \varepsilon ^{-1}Z_t$ as $\varepsilon \rightarrow 0$ in the
situation, where $\delta=\delta _\varepsilon  \rightarrow 0$ but {\it slowly.}
What we need for the
consistency of the estimator $\vartheta^* _{\delta_\varepsilon  ,\varepsilon
}$ is the condition: for any $\nu >0$
$$
\lim_{\varepsilon \rightarrow 0}\varepsilon ^{-2}\inf_{\left|\vartheta
  -\vartheta _0\right|>\nu }\int_{0}^{\delta _\varepsilon
}\left[x_t\left(\vartheta \right)-x_t\left(\vartheta_0 \right)\right]^2{\rm
  d}t\rightarrow \infty .
$$
For example, if the derivative $\left|\dot S\left(\vartheta_0
,0,x_0\right)\right|\geq \gamma >0 $, then for the function
$$
\dot x_t\left(\vartheta
\right)=\int_{0}^{t}\exp\left\{\int_{s}^{t}S'_x\left(\vartheta
,v,x_v\right){\rm d}v\right\} \dot S\left(\vartheta ,s,x_s\right)\,{\rm d}s
$$
and small $\delta _\varepsilon $ we have
$$
\varepsilon ^{-2}\inf_{\left|\vartheta -\vartheta _0\right|>\nu
}\int_{0}^{\delta _\varepsilon }\left[x_t\left(\vartheta
  \right)-x_t\left(\vartheta_0 \right)\right]^2{\rm d}t\geq \frac{\gamma
  ^2\,\delta _\varepsilon ^3\,\nu ^2}{6\,\varepsilon ^{2}} .
$$
Let us consider the linear case
$$
{\rm d}X_t=\vartheta X_t{\rm d}t+\varepsilon {\rm d}W_t,\quad X_0=x_0>0,\quad
0\leq t\leq T.
$$
Then the MLE can be written explicitly
$$
\hat\vartheta _{t,\varepsilon} =\frac{\int_{0}^{t}X_s{\rm
    d}X_s}{\int_{0}^{t}X_s^2{\rm d}s} =\vartheta +\varepsilon \frac{\int_{0}^{t}X_s{\rm
    d}W_s}{\int_{0}^{t}X_s^2{\rm d}s}
$$
and
\begin{align*}
\hat\vartheta _{\delta _\varepsilon ,\varepsilon}-\vartheta =\varepsilon \frac{\int_{0}^{\delta _\varepsilon }X_s{\rm
    d}W_s}{\int_{0}^{\delta _\varepsilon }X_s^2{\rm d}s}\;\sim \;
\frac{\varepsilon \,W_{\delta _\varepsilon }}{x_0\;\delta _\varepsilon }\;\sim
\; \frac{\varepsilon \,W_{1}}{x_0\;\delta _\varepsilon ^{1/2}} .
\end{align*}
Therefore, if $\varepsilon \delta _\varepsilon ^{-1/2}\rightarrow 0$ (for
example, $\delta _\varepsilon =\varepsilon ^2\ln\frac{1}{\varepsilon }$) then
$\hat Y_t\rightarrow Y_t$ for all $t\in \left[\delta _\varepsilon ,T\right]$.

\bigskip

{\it Approximation of the BSDE.}
Note that $\hat Y_t$ is approximation of the solution of the BSDE \eqref{6},
but the stochastic process $\hat Y_t $ itself satisfies another stochastic
differential equation. To simplify the notations let us put
\begin{align*}
{\rm I}_t&={\rm I}\left(\vartheta ^*_{\delta ,\varepsilon },x^t\left(\vartheta
^*_{\delta ,\varepsilon }\right)\right),\qquad \quad
\Delta _t=\Delta _t\left(\vartheta ^*_{\delta ,\varepsilon },x^t_\delta
\right)+\Delta _\delta \left(\vartheta ^*_{\delta ,\varepsilon },x^\delta
\right) ,\\
b_t\left(x\right)&=\frac{\dot S\left(\vartheta ^*_{\delta ,\varepsilon
  },t,x\right)}{\sigma \left(t,x\right)},\qquad \qquad c_t\left(x\right)=\frac{
  S\left(\vartheta_0,t,x\right)   -S\left(\vartheta ^*_{\delta ,\varepsilon
  },t,x\right) }{\sigma \left(t,x\right)}.
\end{align*}

Then we can write the stochastic differential of the  one-step MLE
$\tilde\vartheta _{t,\varepsilon }$
as follows
\begin{align*}
{\rm d}\tilde\vartheta _{t,\varepsilon }= {\rm
  I}_t^{-1}b_t\left(X_t\right) \left[c_t\left(X_t\right)-{\rm
  I}_t^{-1}b_t\left(x_t\right)\Delta _t  \right]  {\rm d}t+
\varepsilon {\rm I}_t^{-1}b_t\left(X_t\right){\rm d}W_t,\quad \tilde\vartheta _{\delta ,\varepsilon },
\end{align*}
where $t\in\left[\delta ,T\right]$.

The approximation of the BSDE is the following equation (below
$u=u\left(t,X_t,\tilde\vartheta _{t,\varepsilon }
\right)$ and $\delta \leq t\leq T$)
\begin{align*}
{\rm d}\hat Y_t&=u'_t\;{\rm d}t+u'_x\;S\left(\vartheta _0,t,X_t\right){\rm d}t+\frac{\varepsilon
  ^2}{2}u''_{x,x}\;\sigma
\left(t,X_t\right)^2 {\rm d}t+\varepsilon u'_x\;\sigma\left(t,X_t\right){\rm d}W_t\\
&\qquad +\dot u\; {\rm
  I}_t^{-1}b_t\left(X_t\right) \left[c_t\left(X_t\right)-{\rm
  I}_t^{-1}b_t\left(x_t\right)\Delta _t  \right]  {\rm d}t +\frac{\varepsilon
  ^2}{2}\ddot u\;  {\rm
  I}_t^{-2}b_t\left(X_t\right)^2   {\rm d}t\\
&\qquad  +\varepsilon\dot
u\; {\rm
  I}_t^{-1}b_t\left(X_t\right){\rm d}W_t +\frac{\varepsilon ^2}{2}\dot
u'_x\;  {\rm
  I}_t^{-1}b_t\left(X_t\right)\sigma\left(t,X_t\right) {\rm d}t,\qquad  \hat Y_\delta  .
\end{align*}
It can be written as follows
\begin{align*}
{\rm d}\hat Y_t&=-f\left(t,X_t,\hat Y_t,\hat Z_t\right){\rm d}t+\hat Z_t{\rm d}W_t\\
&\qquad +u'_x\,S\left(\vartheta _0,t,X_t\right){\rm d}t +\dot u\;  {\rm
  I}_t^{-1}b_t\left(X_t\right) \left[c_t\left(X_t\right)-{\rm
  I}_t^{-1}b_t\left(x_t\right)\Delta _t  \right]  {\rm d}t\\
&\qquad +\frac{\varepsilon ^2}{2}\ddot u\; {\rm
  I}_t^{-2}b_t\left(X_t\right)^2   {\rm d}t  +\varepsilon\dot
u\; {\rm
  I}_t^{-1}b_t\left(X_t\right){\rm d}W_t\\
&\qquad +\frac{\varepsilon ^2}{2}\dot
u'_x\; {\rm
  I}_t^{-1}b_t\left(X_t\right)\sigma\left(t,X_t\right) {\rm d}t,\qquad  \hat
Y_\delta ,\quad \delta \leq t\leq T.
\end{align*}

\bigskip

{\it Other estimators of $Y_t$.}  There are many possibilities to construct
estimators $\bar Y_t$ of the random function $Y_t$ such that $\varepsilon
^{-1}\left(\bar Y_t-Y_t\right)\Rightarrow {\cal N}$. We can put $\bar
Y_t=u\left(t,X_t,\vartheta ^*_{\delta ,\varepsilon }\right),t\in \left[\delta
  ,T\right]$. Then
$$
\varepsilon^{-1}\left(\bar Y_t-Y_t\right)\Longrightarrow {\cal N}\left(0,\dot
u\left(t,x_t,\vartheta _0\right)^2D_\delta \left(\vartheta _0\right)^2\right) ,
$$
where
\begin{align*}
D_\delta \left(\vartheta _0\right)^2=\int_{0}^{\delta }\frac{\sigma
  \left(s,x_s\right)^2}{\psi \left(s,\vartheta _0\right)} \left(
\int_{s}^{\delta }\psi \left(t,\vartheta _0 \right) \,\dot x_t\,{\rm d}t \right)^2{\rm d}s
\end{align*}
and
$$
\psi \left(t,\vartheta _0\right)=\exp\left\{\int_{0}^{t}
S'_x\left(s,x_s\right)\,{\rm d}s\right\}
$$
 (see Theorem 7.5, \cite{Kut94}). It is known that for $t\geq \delta $
$$
D_\delta
\left(\vartheta _0\right)^2\geq  {\rm I}\left(\vartheta _0,x^\delta \left(\vartheta
_0\right)\right)^{-1}\geq  {\rm I}\left(\vartheta _0,x^t\left(\vartheta
_0\right)\right)^{-1}.
$$
Therefore this approximation is not asymptotically efficient.

Another possibility is to use the limit equation ($\varepsilon =0$)
\begin{align*}
\frac{\partial u^{0}}{\partial t}+S\left(\vartheta,x\right) \frac{\partial
  u^{0}}{\partial x}=-f\left(x,u^{0},0\right),\qquad
u^{0}\left(T,x,\vartheta \right)=\Phi \left(x\right)
\end{align*}
and to put $\bar Y_t=u^{0}\left(t,X_t,\vartheta^*_{\delta
  ,\varepsilon } \right)$. Under regularity conditions we have the convergence
$$
\sup_\vartheta \left|u\left(t,X_t,\vartheta
\right)-u\left(t,x_t,\vartheta \right)\right|\rightarrow
0,\qquad u\left(t,x,\vartheta \right)\longrightarrow
u^{0}\left(t,x,\vartheta\right).
$$
Therefore $\bar Y_t-Y_t\rightarrow 0 $.   This means that the both random
functions have the same (deterministic) limit. Therefore such solutions are
not asymptotically efficient and are essentially less interesting.

\bigskip

{\it Linear case. }
Let us consider one example. Suppose that
$$
{\rm d}X_t=\vartheta {\rm d}t+\varepsilon\sigma  {\rm d}W_t,\quad X_0=x_0,
\quad 0\leq t\leq T,
$$
where $\vartheta \in\Theta =\left(a ,b \right)$ and we are given two
functions $f\left(x,z\right)=\beta y+\gamma z$ and  $\Phi
\left(x\right)$. The variables $\sigma ,\beta , \gamma $ are known
constants and $\vartheta $ is unknown parameter. The function $\Phi
\left(x\right)$ has two continuous derivatives with polynomial
majorants.    We have to construct the  BSDE
\begin{equation}
\label{18}
{\rm d}Y_t=-\left(\beta Y_t+\gamma Z_t\right){\rm d}t+Z_t{\rm d}W_t,\qquad
Y_T=\Phi \left(X_T\right).
\end{equation}
The corresponding PDE is
 \begin{equation}
\label{19}
  \left\{\begin{aligned} &\frac{\partial u}{\partial
   t}+\frac{1}{2}\varepsilon^2\sigma^2\frac{\partial^2 u}{\partial
   x^2} +(\vartheta+\varepsilon\sigma\gamma)\frac{\partial
       u}{\partial x}+\beta u=0,\ 0\leq t\leq T,
    \\
&u(T,x,\vartheta )=\Phi(x),\ x\in\mathbb{R}. \\
      \end{aligned}\right.
 \end{equation}
  with the solution
\begin{align*}
  u(t,x,\vartheta)= e^{\beta(T-t)} \; G(t,x,\vartheta),
\end{align*}
where we denoted
\begin{align*}
 G(t,x,\vartheta)=\int_{-\infty}^\infty e^{-\frac{z^2
   }{2\varepsilon^2\sigma^2(T-t)}}
 \; \frac{\Phi(x+(\vartheta+\varepsilon\sigma\gamma)(T-t)-z)}{\sqrt{2\pi
     \varepsilon^2\sigma^2(T-t)}    }\: {\rm d} z
\end{align*}
 Then we can put
        \begin{eqnarray*}
   &&Y_t=u(t,X_t,\vartheta)=e^{\beta(T-t)}\;G(t,X_t,\vartheta),\\
     &&Z_t=\varepsilon\;\sigma\;  u'(t,X_t,\vartheta)=\varepsilon\;\sigma\;
     e^{\beta(T-t)}\;G_x'(t,X_t,\vartheta),
        \end{eqnarray*}
  and obtain  the BSDE \eqref{18}.

  Note that
 \begin{eqnarray*}
    &&G'_x(t,x,\vartheta)=\int_{-\infty}^\infty       e^{-\frac{z^2}{2\varepsilon
     ^2\sigma^2(T-t)}}\;\frac{\Phi'(x+(\vartheta+\varepsilon\sigma\gamma)(T-t)-z)}{\sqrt{2\pi
   \varepsilon^2\sigma^2(T-t)}               }\;{\rm d} z.
 \end{eqnarray*}
   We have
   \begin{align*}
\dot u_\vartheta(t,x,\vartheta)&=(T-t)e^{\beta(T-t)}   \int_{-\infty}^\infty
   e^{-\frac{z^2}{2\sigma^2(T-t)}}\;\frac{\Phi'(x+(\vartheta+\varepsilon\sigma\gamma)(T-t)-z)}{
    \sqrt{2\pi\varepsilon^2\sigma^2(T-t)}   }\;{\rm d}z\\
  & =(T-t)e^{\beta(T-t)}G'_x(t,x,\vartheta),
\end{align*}
   and
  \begin{eqnarray*}
  &&\dot u'(t,x,\vartheta)    =(T-t)e^{\beta(T-t)}G''(t,x,\vartheta),\\
  &&\ddot u(t,x,\vartheta)      =(T-t)^2e^{\beta(T-t)}G''(t,x,\vartheta).
        \end{eqnarray*}
In this model the MLE $\hat\vartheta _{t,\varepsilon }$ can be explicitly
written
$$
\hat\vartheta _{t,\varepsilon }=\frac{X_t}{t} =\vartheta _0+\varepsilon \sigma
\frac{W_t}{t}\sim {\cal N}\left(\vartheta _0,\frac{\varepsilon ^2\sigma ^2}{t}\right)
$$
and for all $t\in(0,T]$ is consistent. Therefore we can put
\begin{align*}
\hat Y_t&=e^{\beta(T-t)}\;G(t,X_t,\hat\vartheta _{t,\varepsilon }),\qquad t\in(0,T]\\
\hat Z_t&=\varepsilon\;\sigma\; e^{\beta(T-t)}\;G_x'(t,X_t,\hat\vartheta
_{t,\varepsilon }),\qquad t\in(0,T]
\end{align*}
and according to the Theorem \ref{T3} this approximation is asymptotically
efficient.

The limit solution of PDE is
$$
u^0\left(t,x,\vartheta \right)=e^{\beta \left(T-t\right)}\,\Phi
\left(x+\vartheta \left(T-t\right)\right).
$$
Hence
\begin{align*}
\lim_{\varepsilon \rightarrow 0}\varepsilon ^{-2}\Ex_\vartheta \left(\hat
Y_t-Y_t\right)^2=\left(T-t\right)^2t ^{-1}\,\sigma ^{2}e^{2\beta \left(T-t\right)} \Phi '\left(x_0+\vartheta T\right)^2.
\end{align*}
The last expression is also the rhs in the lower bound  \eqref{15}.

\bigskip

{\it On regularity condition \eqref{q}.} The most difficult to verify are the
conditions imposed on $u\left(t,x,\vartheta \right)$ (regularity
w.r.t. $\vartheta $, convergence to $u^0\left(t,x,\vartheta \right)$ and
majorations in $x$ of $u\left(\cdot \right)$ and its derivatives).  Sufficient conditions for the convergence of the solution
$u\left(t,x,\vartheta \right)$ in homogeneous case and a particular case of
the linear $f\left(\cdot \right)$ are given in the  Theorem 1.3.1 in
\cite{FW}.

\bigskip

{\it Some generalizations.} The parameter $\vartheta $ in our work
is supposed to be one-dimensional and the loss function $\ell
\left(y\right)$ in the theorems \ref{T2} and \ref{T3} quadratic,
$\ell \left(y\right)=y^2$. The case of multidimensional parameter
and more general class of loss functions can be treated by a similar
way, but in this case we have to use Hajek-Le Cam lower bound as
follows. Suppose that $\vartheta \in \Theta $ where $\Theta $ is an
open bounded subset of $\RR^d$ and the corresponding regularity
conditions are fulfilled. Then the family of measures
$\left\{\Pb_\vartheta ^{\left(t\right)},\vartheta \in \Theta
\right\}$ (induced in the space of realizations by the stochastic
processes $X^t=\left(X_s,0\leq s\leq t\right)$) is {\it locally
asymptotically normal} (LAN), i.e., the normalized likelihood ratio
function $Z_{t,\varepsilon} \left(v\right)=L\left(\vartheta
_0+\varepsilon v,\vartheta _0, X^t\right)$ admits the representation
\begin{align*}
Z_\varepsilon \left(v\right)=\exp\left\{\langle v,\Delta_t\left(\vartheta
_0,X_0^t\right)\rangle -\frac{1}{2} v^*{\rm I}\left(\vartheta
_0,x^t\left(\vartheta _0\right)\right)v+r_\varepsilon  \right\} ,\qquad v\in \RR^d,
\end{align*}
where $\Delta_t\left(\vartheta_0,X_0^t\right)\Rightarrow
{\cal N}\left(0, {\rm I}\left(\vartheta_0,x^t\left(\vartheta _0\right)\right)\right)$
and $r_\varepsilon \rightarrow 0$. Suppose that the loss function
$\ell\left(y\right)$ is nonnegative, continuous at point 0 and
$\ell\left(0\right)=0$, but is not identically 0, is symmetric
$\ell\left(y\right)=\ell\left(\left|y\right|\right)$ and
$\ell\left(y\right),y\geq 0$  is nondecreasing. Then we have the
following Hajek-Le Cam lower bound: for all $\vartheta _0\in \Theta $ and for
all estimators $\bar Y_t, t\in \left[\delta ,T\right]$
\begin{align*}
&\Liminf_{\nu \rightarrow 0}\Liminf_{\varepsilon \rightarrow 0} \sup_{\left|\vartheta
  -\vartheta _0\right|\leq \nu } \Ex_\vartheta \ell\left(\varepsilon ^{-1}\left( \bar
Y_t-Y_t\right)\right)\geq \Ex_{\vartheta_0} \ell\left(\langle\dot u^0\left(t,x_t\left(\vartheta
  _0\right),\vartheta _0\right),\xi _t\left(\vartheta _0\right)\rangle\right) .
\end{align*}

For the proof of LAN see Lemma 2.2 in  \cite{Kut94} and for the lower bound see
Theorem 2.12.1 in  \cite{IH81} (we have to modify slightly the proof, because
we estimate not $\vartheta $ but some random function of $\vartheta $).

The next step, of course, is to prove the asymptotic efficiency of the
one-step MLE in the sense of this bound.

\end{document}